\documentclass[12pt]{article}

\makeatletter

\@addtoreset{equation}{section}
\makeatother

\newtheorem{definition}{Definition}[section]
\newtheorem{theorem}{Theorem}[section]
\newtheorem{proposition}{Proposition}[section]
\newtheorem{corollary}{Corollary}[section]
\newtheorem{lemma}{Lemma}[section]
\newtheorem{remark}{Remark}[section]

\usepackage{amssymb}
\usepackage{latexsym}

\begin{document} 

\title{Generalized Uncertainty Relation Associated with a Monotone or an Anti-Monotone Pair Skew Information}
\author{Kenjiro Yanagi\thanks{
Division of Applied Mathematical Science,
Graduate School of Science and Engineering,
Yamaguchi University,
Ube, 755-8611 Japan. 
E-mail: {\tt yanagi@\allowbreak
yamaguchi-u.\allowbreak
ac.\allowbreak
jp},  
This research was partially supported by the Ministry of Education, Science, Sports and Culture, Grant-in-Aid for 
Scientific Research (C), 23540208} 
and 
Satoshi Kajihara\thanks{
Graduate School of Science and Engineering,
Yamaguchi University,
Yamaguchi, 753-8511 Japan.
E-mail: {\tt p002va@\allowbreak
yamaguchi-u.\allowbreak
ac.\allowbreak
jp}}}
\date{}

\maketitle

\begin{flushleft}
{\bf Abstract.} We give a trace inequality related to the uncertainty relation 
based on the monotone or anti-monotone pair skew information which is one of generalizations of result given by 
\cite{KoYo:monotone}. And it includes the result for generalized Wigner-Yanase-Dyson skew information as a particular case 
(\cite{Ya:GWYD1}).
\end{flushleft}

\begin{flushleft}
{\bf Key Words:} Uncertainty relation, Wigner-Yanase-Dyson skew information
\end{flushleft}

\section{Introduction}

Wigner-Yanase skew information
\begin{eqnarray*}
I_{\rho}(H) & = & \frac{1}{2}Tr \left[ \left( i \left[ \rho^{1/2},H \right] \right)^2 \right] \\
& = & Tr[\rho H^2]-Tr[\rho^{1/2}H\rho^{1/2}H] 
\end{eqnarray*}
was defined in \cite{WY:inf}. This quantity can be considered as a kind of 
the degree for non-commutativity between a quantum state $\rho$ and an 
observable $H$. Here we denote the commutator by $[X,Y] = XY-YX$. This 
quantity was generalized by Dyson 
\begin{eqnarray*}
I_{\rho,\alpha}(H) & = & \frac{1}{2}Tr[(i[\rho^{\alpha},H])(i[\rho^{1-\alpha},H])] \\
& = & Tr[\rho H^2]-Tr[\rho^{\alpha}H\rho^{1-\alpha}H],  \alpha \in [0,1]
\end{eqnarray*}
which is known as the Wigner-Yanase-Dyson skew information. It is famous that 
the convexity of $I_{\rho,\alpha}(H)$ with respect to $\rho$ was successfully 
proven by E.H.Lieb in \cite{Li:convex}. 
And also this quantity was generalized by Cai and Luo
\begin{eqnarray*}
&   & I_{\rho,\alpha,\beta}(H) \\
& = & \frac{1}{2}Tr[(i[\rho^{\alpha},H])(i[\rho^{\beta},H])\rho^{1-\alpha-\beta}] \\
& = & \frac{1}{2} \{ Tr[\rho H^2]+Tr[\rho^{\alpha+\beta}H\rho^{1-\alpha-\beta}H]-Tr[\rho^{\alpha}H\rho^{1-\alpha}H]-Tr[\rho^{\beta}H\rho^{1-\beta}H] \},
\end{eqnarray*}
where $\alpha, \beta \geq 0, \alpha+\beta \leq 1$. The convexity of $I_{\rho,\alpha,\beta}(H)$ with respect to $\rho$ 
was proven by Cai and Luo in \cite{CaLu:convex} under some restrictive condition. 
In this paper we let $M_n(\mathbb{C})$ be the set of all $n \times n$ complex matrices,  
$M_{n,sa}(\mathbb{C})$ be the set of all $n \times n$ self-adjoint matrices, 
$M_{n,+}(\mathbb{C})$ be the set of strictly positive elements of $M_n(\mathbb{C})$ and 
$M_{n,+,1}(\mathbb{C})$ be the set of strictly positive density matrices, that is 
$M_{n,+,1}(\mathbb{C}) = \{ \rho \in M_n(\mathbb{C}) | Tr[\rho] = 1, \rho > 0 \}$. 
If it is not otherwise specified, from now on we shall treat the case of faithful states, 
that is $\rho > 0$. The relation between the Wigner-Yanase skew 
information and the uncertainty relation was studied in \cite{LuZh:skew}. 
Moreover the relation between the Wigner-Yanase-Dyson skew information and the 
uncertainty relation was studied in \cite{Ko:matrix, YaFuKu:gen}. In our paper 
\cite{YaFuKu:gen} and \cite{Ya:WYD}, we defined a generalized skew information and then derived 
a kind of an uncertainty relations. And also in \cite{Ya:GWYD1} and \cite{Ya:GWYD2} , 
we gave an uncertainty relation of two parameter generalized Wigner-Yanase-Dyson skew information. 
In this paper, we consider three parameter generalized Wigner-Yanase-Dyson skew information and 
give a kind of generalized uncertainty relations which is a generalization of the result of Ko and Yoo \cite{KoYo:monotone}.

\section{Trace inequality of Wigner-Yanase-Dyson skew information}

We review the relation between the Wigner-Yanase skew information and the uncertainty relation. 
In quantum mechanical system, the expectation value of an observable $H$ 
in a quantum state $\rho$ is expressed by $Tr[\rho H]$. It is natural 
that the variance for a quantum state $\rho$ and an observable $H$ is defined by 
$V_{\rho}(H) = Tr[\rho (H-Tr[\rho H]I)^2] = Tr[\rho H^2]-Tr[\rho H]^2$. 
It is famous that we have 
\begin{equation}
V_{\rho}(A)V_{\rho}(B) \geq \frac{1}{4}|Tr[\rho [A,B]]|^2 
\label{eq:num2-1}
\end{equation}
for a quantum state $\rho$ and two observables $A$ and $B$. The further strong results 
was given by Schr\"odinger 
$$
V_{\rho}(A)V_{\rho}(B)-|Re\{Cov_{\rho}(A,B)\}|^2 \geq \frac{1}{4}|Tr[\rho [A,B]]|^2, 
$$
where the covariance is defined by 
$Cov_{\rho}(A,B) = Tr[\rho (A-Tr[\rho A]I)(B-Tr[\rho B]I)]$. 
However, the uncertainty relation for the Wigner-Yanase skew information failed. 
(See \cite{LuZh:skew, Ko:matrix, YaFuKu:gen}) 
$$
I_{\rho}(A) I_{\rho}(B) \geq \frac{1}{4}|Tr[\rho [A,B]]|^2.
$$
Recently, S.Luo introduced the quantity $U_{\rho}(H)$ representing a quantum uncertainty 
excluding the classical mixture: 
\begin{equation}
U_{\rho}(H) = \sqrt{V_{\rho}(H)^2-(V_{\rho}(H)-I_{\rho}(H))^2}, 
\label{eq:num2-2}
\end{equation}
then he derived the uncertainty relation on $U_{\rho}(H)$ in \cite{Lu:Hei}: 
\begin{equation}
U_{\rho}(A) U_{\rho}(B) \geq \frac{1}{4}|Tr[\rho [A,B]]|^2.
\label{eq:num2-3}
\end{equation}
Note that we have the following relation 
\begin{equation}
0 \leq I_{\rho}(H) \leq U_{\rho}(H) \leq V_{\rho}(H). 
\label{eq:num2-4}
\end{equation}
The inequality (\ref{eq:num2-3}) is a refinement of the inequality (\ref{eq:num2-1}) 
in the sense of (\ref{eq:num2-4}). 
In \cite{Ya:WYD}, we studied one-parameter extended inequality for the inequality (\ref{eq:num2-3}). 

\begin{definition}
For $0 \leq \alpha \leq 1$, a quantum state $\rho$ and an observable $H$, we define 
the Wigner-Yanase-Dyson skew information 
\begin{eqnarray*}
I_{\rho,\alpha}(H) & = & \frac{1}{2}Tr[(i [\rho^{\alpha},H_0])(i [\rho^{1-\alpha},H_0])] \\
& = & Tr[\rho H_0^2]-Tr[\rho^{\alpha}H_0\rho^{1-\alpha}H_0] 
\end{eqnarray*}
and we also define 
\begin{eqnarray*}
J_{\rho,\alpha}(H) & = & \frac{1}{2}Tr[\{ \rho^{\alpha},H_0 \} \{ \rho^{1-\alpha},H_0 \}] \\
& = & Tr[\rho H_0^2]+Tr[\rho^{\alpha}H_0 \rho^{1-\alpha}H_0], 
\end{eqnarray*}
where $H_0 = H-Tr[\rho H]I$ and we denote the anti-commutator by $\{ X,Y \} = XY+YX$.  
\label{def:definition2-1}
\end{definition}

Note that we have 
$$
\frac{1}{2}Tr[(i [\rho^{\alpha},H_0])(i [\rho^{1-\alpha},H_0])] = \frac{1}{2}Tr[(i [\rho^{\alpha},H])(i [\rho^{1-\alpha},H])]
$$
but we have 
$$
\frac{1}{2}Tr[\{ \rho^{\alpha},H_0 \} \{ \rho^{1-\alpha},H_0 \}] \neq \frac{1}{2}Tr[\{\rho^{\alpha},H \} \{ \rho^{1-\alpha},H \}]. 
$$
Then we have the following inequalities: 
\begin{equation}
I_{\rho,\alpha}(H) \leq I_{\rho}(H) \leq J_{\rho}(H) \leq J_{\rho,\alpha}(H), 
\label{eq:num2-5}
\end{equation}
since we have $Tr[\rho^{1/2}H\rho^{1/2}H] \leq Tr[\rho^{\alpha}H\rho^{1-\alpha}H]$. 
(See \cite{Bo:some,Fu:trace} for example.) If we define 
\begin{equation}
U_{\rho,\alpha}(H) = \sqrt{V_{\rho}(H)^2-(V_{\rho}(H)-I_{\rho,\alpha}(H))^2}, 
\label{eq:num2-6}
\end{equation}
as a direct generalization of Eq.(\ref{eq:num2-2}), then we have 
\begin{equation}
0 \leq I_{\rho,\alpha}(H) \leq U_{\rho,\alpha}(H) \leq U_{\rho}(H) 
\label{eq:num2-7}
\end{equation}
due to the first inequality of (\ref{eq:num2-5}). We also have 
$$
U_{\rho,\alpha}(H) = \sqrt{I_{\rho,\alpha}(H) J_{\rho,\alpha}(H)}. 
$$
From the inequalities (\ref{eq:num2-4}),(\ref{eq:num2-6}),(\ref{eq:num2-7}), our situation is that we have 
$$
0 \leq I_{\rho,\alpha}(H) \leq I_{\rho}(H) \leq U_{\rho}(H)
$$
and
$$
0 \leq I_{\rho,\alpha}(H) \leq U_{\rho,\alpha}(H) \leq U_{\rho}(H).
$$
We gave the following uncertainty relation with respect to $U_{\rho,\alpha}(H)$ as a direct 
generalization of the inequality (\ref{eq:num2-3}). 

\begin{theorem}[\cite{Ya:WYD}]
For $0 \leq \alpha \leq 1$, a quantum state $\rho$ and observables $A,B$, 
\begin{equation}
U_{\rho,\alpha}(A) U_{\rho,\alpha}(B) \geq \alpha(1-\alpha)|Tr[\rho [A,B]]|^2.
\label{eq:num2-8}
\end{equation}
\label{th:theorem2-1}
\end{theorem}

Now we define the two parameter extensions of Wigner-Yanase skew information and 
give an uncertainty relation under some conditions.  

\begin{definition}
For $\alpha, \beta \geq 0$, a quantum state $\rho$ and an observable $H$, 
we define the generalized Wigner-Yanase-Dyson skew information
\begin{eqnarray*}
&   & I_{\rho,\alpha,\beta}(H) \\
& = & \frac{1}{2}Tr\left[(i[\rho^{\alpha},H_0])(i[\rho^{\beta},H_0])\rho^{1-\alpha-\beta} \right] \\
& = & \frac{1}{2}\{ Tr[\rho H_0^2]+Tr[\rho^{\alpha+\beta}H_0 \rho^{1-\alpha-\beta}H_0]-Tr[\rho^{\alpha}H_0 \rho^{1-\alpha}H_0]-Tr[\rho^{\beta}H_0 \rho^{1-\beta}H_0] \} 
\end{eqnarray*}
and we define 
\begin{eqnarray*}
&   & J_{\rho,\alpha,\beta}(H) \\
& = & \frac{1}{2}Tr \left[\{\rho^{\alpha},H_0 \} \{\rho^{\beta},H_0 \} \rho^{1-\alpha-\beta} \right] \\
& = & \frac{1}{2}\{ Tr[\rho H_0^2]+Tr[\rho^{\alpha+\beta}H_0 \rho^{1-\alpha-\beta}H_0]+Tr[\rho^{\alpha}H_0 \rho^{1-\alpha}H_0]+Tr[\rho^{\beta}H_0 \rho^{1-\beta}H_0] \},
\end{eqnarray*}
where $H_0 = H-Tr[\rho H]I$ and we denote the anti-commutator by $\{X,Y \} = XY+YX$. We remark that $\alpha+\beta = 1$ implies 
$I_{\rho,\alpha}(H) = I_{\rho,\alpha,1-\alpha}(H)$ and $J_{\rho,\alpha}(H) = J_{\rho,\alpha,1-\alpha}(H)$.
We also define 
$$
U_{\rho,\alpha,\beta}(H) = \sqrt{I_{\rho,\alpha,\beta}(H) J_{\rho,\alpha,\beta}(H)}.
$$
\label{def:definition2-2}
\end{definition}

In this paper we assume that $\alpha, \beta \geq 0$ do not necessarily satisfy the condition $\alpha+ \beta \leq 1$. 
We give the following theorem.
 
\begin{theorem}[\cite{Ya:GWYD1}]
For $\alpha, \beta \geq 0$ and $\alpha + \beta \geq 1$ or $\alpha + \beta \leq \frac{1}{2}$ and observables $A, B$, 
\begin{equation}
U_{\rho,\alpha, \beta}(A)U_{\rho,\alpha,\beta}(B) \geq \alpha \beta|Tr[\rho [A,B]]|^2.
\label{eq:num2-9}
\end{equation}
\label{th:theorem2-2}
\end{theorem}

And we also define the two parameter extensions of Wigner-Yanase skew information which are different from 
Definition \ref{def:definition2-2}.

\begin{definition}
For $\alpha, \beta \geq 0$, a quantum state $\rho$ and an observable $H$, we define the generalized 
Wigner-Yanase-Dyson skew information 
\begin{eqnarray*}
&   & \tilde{I}_{\rho,\alpha,\beta}(H) \\
& = & \frac{1}{2}Tr \left[(i[\rho^{\alpha},H_0])(i[\rho^{\beta},H_0]) \right] \\
& = & Tr[\rho^{\alpha+\beta}H_0^2]-Tr[\rho^{\alpha}H_0 \rho^{\beta}H_0]. 
\end{eqnarray*}
and we define 
\begin{eqnarray*}
&   & \tilde{J}_{\rho,\alpha,\beta}(H) \\
& = & \frac{1}{2}Tr \left[ \{ \rho^{\alpha},H_0 \} \{ \rho^{\beta},H_0 \} \right] \\
& = & Tr[\rho^{\alpha+\beta}H_0^2]+Tr[\rho^{\alpha}H_0 \rho^{\beta}H_0], 
\end{eqnarray*}
where $H_0 = H-Tr[\rho H]I$ and we denote the anti-commutator by $\{ X,Y \} = XY+YX$. We remark that 
$\alpha+\beta = 1$ implies $I_{\rho,\alpha}(H) = \tilde{I}_{\rho,\alpha,1-\alpha}(H)$ and 
$J_{\rho,\alpha}(H) = \tilde{J}_{\rho,\alpha,1-\alpha}(H)$. We also define 
$$
\tilde{U}_{\rho,\alpha,\beta}(H) = \sqrt{\tilde{I}_{\rho,\alpha,\beta}(H) \tilde{J}_{\rho,\alpha,\beta}(H)}.
$$
\label{def:definition2-3}
\end{definition}

Then we give the following theorem. 

\begin{theorem}[\cite{Ya:GWYD2}]
For $\alpha, \beta \geq 0~(\alpha \beta \neq 0)$ and observables $A, B$, 
$$
\tilde{U}_{\rho,\alpha,\beta}(A) \tilde{U}_{\rho,\alpha,\beta}(B) \geq \frac{\alpha \beta}{(\alpha+\beta)^2}|Tr[\rho^{\alpha+\beta}[A,B]]|^2. 
$$
\label{th:theorem2-3}
\end{theorem}

\begin{remark}
We remark that (\ref{eq:num2-8}) is derived by putting $\beta = 1-\alpha$ in (\ref{eq:num2-9}). Then Theorem \ref{th:theorem2-2} 
is a generalization of Theorem \ref{th:theorem2-1} given in \cite{Ya:WYD}. 
\label{re:remark2-1}
\end{remark}

\section{Trace inequality of monotone or anti-monotone pair skew information}

\begin{definition}
Let $f(x), g(x)$ be nonnegative continuous functions defined on the interval $[0,1]$.
We call the pair $(f,g)$ a compatible in log-increase, monotone pair (CLI monotone pair, in short) if
\begin{description}
\item[(a)] $(f(x)-f(y))(g(x)-g(y)) \geq 0$ for all $x, y \in [0,1]$. 
\item[(b)] $f(x), g(x)$ are differentiable on $(0,1)$ and 
$$
0 \leq \inf_{0<x<1} \frac{G^{'}(x)}{F^{'}(x)} \leq \sup_{0<x<1} \frac{G^{'}(x)}{F^{'}(x)} < \infty,
$$
where $F(x) = \log f(x), G(x) = \log g(x)$.
\end{description}
\label{def:definition3-1}
\end{definition}

\begin{definition}
Let $f(x), g(x)$ be nonnegative continuous functions defined on the interval $[0,1]$.
We call the pair $(f,g)$ a compatible in log-increase, anti-monotone pair (CLI anti-monotone pair, in short) if
\begin{description}
\item[(a)] $(f(x)-f(y))(g(x)-g(y)) \leq 0$ for all $x, y \in [0,1]$. 
\item[(b)] $f(x), g(x)$ are differentiable on $(0,1)$ and 
$$
-\infty < \inf_{0<x<1} \frac{G^{'}(x)}{F^{'}(x)} \leq \sup_{0<x<1} \frac{G^{'}(x)}{F^{'}(x)} \leq 0,
$$
where $F(x) = \log f(x), G(x) = \log g(x)$.
\end{description}
\label{def:definition3-2}
\end{definition}

Let $f(x), g(x), h(x)$ be nonnegative continuous functions defined on $[0,1]$ and be differentiable on $(0,1)$.
We assume that $(f,g)$ is CLI monotone pair and $(f,h)$ is CLI monotone or anti-monotone pair. 
We introduce the correlation functions in the following way. 

\begin{definition}
\begin{eqnarray*}
I_{\rho,(f,g,h)}(H) & = & \frac{1}{2}Tr[(i[f(\rho),H_0])(i[g(\rho),H_0])h(\rho)] \\
& = & -\frac{1}{2}Tr[(f(\rho),H_0])([g(\rho),H_0])h(\rho)] \\
& = & -\frac{1}{2}Tr[(f(\rho)H_0-H_0f(\rho))(g(\rho)H_0-H_0g(\rho))h(\rho)] \\
& = & -\frac{1}{2}Tr[f(\rho)H_0g(\rho)H_0h(\rho)-f(\rho)H_0^2g(\rho)h(\rho)] \\
&   & +\frac{1}{2}Tr[H_0f(\rho)g(\rho)H_0h(\rho)-H_0f(\rho)H_0g(\rho)h(\rho)] \\
& = & -\frac{1}{2}Tr[f(\rho)h(\rho)H_0g(\rho)H_0-f(\rho)g(\rho)h(\rho)H_0^2] \\
&   & +\frac{1}{2}Tr[f(\rho)g(\rho)H_0h(\rho)H_0-g(\rho)h(\rho)H_0f(\rho)H_0] \\
& = & \frac{1}{2} \{ Tr[f(\rho)g(\rho)h(\rho)H_0^2]+Tr[f(\rho)g(\rho)H_0h(\rho)H_0] \} \\
&   & -\frac{1}{2} \{ Tr[f(\rho)H_0g(\rho)h(\rho)H_0]+Tr[g(\rho)H_0f(\rho)h(\rho)H_0] \}.
\end{eqnarray*}
\begin{eqnarray*}
J_{\rho,(f,g,h)}(H) & = & \frac{1}{2}Tr[ \{f(\rho),H_0 \} \{ g(\rho),H_0 \} h(\rho)] \\
& = & \frac{1}{2}Tr[(f(\rho)H_0+h_0f(\rho))(g(\rho)H_0+H_0g(\rho))h(\rho)] \\
& = & \frac{1}{2}Tr[f(\rho)H_0g(\rho)H_0h(\rho)+f(\rho)H_0^2g(\rho)h(\rho)] \\
&   & +\frac{1}{2}Tr[H_0f(\rho)g(\rho)H_0h(\rho)+H_0f(\rho)H_0g(\rho)h(\rho)] \\
& = & \frac{1}{2} \{ Tr[f(\rho)g(\rho)h(\rho)H_0^2]+Tr[f(\rho)g(\rho)H_0h(\rho)H_0] \} \\
&   & + \frac{1}{2} \{ Tr[f(\rho)H_0g(\rho)h(\rho)H_0]+Tr[g(\rho)H_0f(\rho)h(\rho)H_0] \}.
\end{eqnarray*}
\begin{flushleft}
$~~~~~~~U_{\rho,(f,g,h)}(H) = \sqrt{I_{\rho,(f,g,h)}(H) J_{\rho,(f,g,h)}(H)}$.
\end{flushleft}
\label{def:definition3-3}
\end{definition}

We are ready to state our main result. For $f, g, h$ we let 
\begin{eqnarray}
& & \beta(f,g,h) \nonumber \\
&=& \min \{ \frac{m}{(1+m+n)^2}, \frac{m}{(1+m+N)^2}, \frac{M}{(1+M+n)^2}, \frac{M}{(1+M+N)^2} \}, 
\label{eq:num3-1}
\end{eqnarray}
where
$$
m = \inf_{0<x<1} \frac{G^{'}(x)}{F^{'}(x)}, \; \; \; M = \sup_{0<x<1} \frac{G^{'}(x)}{F^{'}(x)} 
$$
$$
n = \inf_{0<x<1} \frac{H^{'}(x)}{F^{'}(x)}, \; \; \; N = \sup_{0<x<1} \frac{H^{'}(x)}{F^{'}(x)}.
$$
We consider the following two assumptions.
\begin{description}
\item[(I)] $(f,g), (f,h)$ are CLI monotone pair satisfying
$$
1+\frac{G(y)-G(x)}{F(y)-F(x)} \leq \frac{H(y)-H(x)}{F(y)-F(x)}~~{\rm for}~x < y,
$$
where $F(x) = \log f(x), \; G(x) = \log g(x), \; H(x) = \log h(x)$
\item[(II)] $(f,g)$ is CLI monotone pair and $(f,h)$ is CLI anti-monotone pair satisfying 
$$
1+\frac{G(y)-G(x)}{F(y)-F(x)}+\frac{H(y)-H(x)}{F(y)-F(x)} \geq 0~~{\rm for}~x < y.
$$
\end{description}

\begin{theorem}
Under the assumption (I) or (II), the following inequality holds: 
$$
U_{\rho,(f,g,h)}(A) U_{\rho,(f,g,h)}(B) \geq \beta(f,g,h)|Tr[f(\rho)g(\rho)h(\rho)[A,B]] |^2
$$
for $A, B \in M_{n,sa}(\mathbb{C})$.
\label{th:theorem3-1}
\end{theorem}

\section{Proof of Theorem \ref{th:theorem3-1}}

Let $\rho = \sum_{i=1}^n \lambda_i|\phi_i\rangle \langle\phi_i| \in M_{n,+,1}(\mathbb{C})$, where $\{ |\phi_i\rangle \}_{i=1}^n$ 
is an orthonormal set in $\mathbb{C}^n$. Let $(f,g)$ be a CLI monotone pair and $(f,h)$ be a CLI monotone or anti-monotone 
pair. By a simple calculation, we have for any $H \in M_{n,sa}(\mathbb{C})$
\begin{equation}
Tr[f(\rho)g(\rho)h(\rho)H_0^2]  =  \sum_{i,j}\frac{1}{2}\{f(\lambda_i)g(\lambda_i)h(\lambda_i)+f(\lambda_j)g(\lambda_j)h(\lambda_j) \}|a_{ij}|^2. \label{eq:num4-1}
\end{equation}
\begin{equation}
Tr[f(\rho)g(\rho)H_0h(\rho)H_0] = \sum_{i,j}\frac{1}{2}\{f(\lambda_i)g(\lambda_i)h(\lambda_j)+f(\lambda_j)g(\lambda_j)h(\lambda_i) \}|a_{ij}|^2. \label{eq:num4-2} 
\end{equation}
\begin{equation}
Tr[f(\rho)H_0g(\rho)h(\rho)H_0] = \sum_{i,j}\frac{1}{2}\{f(\lambda_i)g(\lambda_j)h(\lambda_j)+f(\lambda_j)g(\lambda_i)h(\lambda_i) \}|a_{ij}|^2. \label{eq:num4-3} 
\end{equation}
\begin{equation}
Tr[g(\rho)H_0f(\rho)h(\rho)H_0] = \sum_{i,j}\frac{1}{2}\{g(\lambda_i)f(\lambda_j)h(\lambda_j)+g(\lambda_j)f(\lambda_i)h(\lambda_i) \}|a_{ij}|^2, \label{eq:num4-4}
\end{equation}
where $a_{ij} = \langle \phi_i|H_0|\phi_j\rangle$ and $a_{ij} = \overline{a_{ji}}$. 
From (\ref{eq:num4-1}) - (\ref{eq:num4-4}), we get
$$
I_{\rho,(f,g,h)}(H) = \frac{1}{2}\sum_{i<j}(f(\lambda_i)-f(\lambda_j))(g(\lambda_i)-g(\lambda_j))(h(\lambda_i)+h(\lambda_j))|a_{ij}|^2.
$$
$$
J_{\rho,(f,g,h)}(H) \geq \frac{1}{2}\sum_{i<j}(f(\lambda_i)+f(\lambda_j))(g(\lambda_i)+g(\lambda_j))(h(\lambda_i)+h(\lambda_j))|a_{ij}|^2.
$$
To prove Theorem \ref{th:theorem3-1}, we need to control a lower bound of a functional coming from a CLI monotone or anti-monotone pair. 
For $f, g, h$ satisfying (I) or (II), we define a function $L$ on $[0,1] \times [0,1]$ by 
\begin{equation}
L(x,y) = \frac{(f(x)^2-f(y)^2)(g(x)^2-g(y)^2)(h(x)+h(y))^2}{(f(x)g(x)h(x)-f(y)g(y)h(y))^2}.
\label{eq:num4-5}
\end{equation}

\begin{proposition}
Under the assumption (I) or (II) 
$$
\min_{x,y \in [0,1]} L(x,y) \geq 16 \beta(f,g,h),
$$
where $\beta(f,g,h)$ is defined in (\ref{eq:num3-1}).
\label{prop:proposition4-1}
\end{proposition}

For the proof of Proposition \ref{prop:proposition4-1}, we need the following lemma. 

\begin{lemma}
If $a, b, c \geq 0$ satisfy $0 < a+b \leq c$ or if $a, b \geq 0, c \leq 0$ satisfy 
$a+b+c > 0$, then the inequality
$$
\frac{(e^{2ar}-1)(e^{2br}-1)(e^{cr}+1)^2}{(e^{(a+b+c)r}-1)^2} \geq \frac{16ab}{(a+b+c)^2}
$$
holds for any real number $r$.
\label{lem:lemma4-1}
\end{lemma}

\begin{flushleft}
{\bf Proof.}  We put $e^r = t$. Then we may prove the following;  
\begin{equation}
(t^{2a}-1)(t^{2b}-1)(t^c+1)^2 \geq \frac{16ab}{(a+b+c)^2}(t^{a+b+c}-1)^2 
\label{eq:num4-6}
\end{equation}
for $t > 0$. It is sufficient to prove (\ref{eq:num4-6}) for $t \geq 1$ and $a, b, c \geq 0, 0 < a+b \leq c$ 
or $a, b \geq 0, c \leq 0, a+b+c > 0$.
\end{flushleft}
By Lemma 3.3 in \cite{Ya:WYD} we have for $0 \leq p \leq 1$ and $s \geq 1$, 
$$
(s^{2p}-1)(s^{2(1-p)}-1) \geq 4p(1-p)(s-1)^2.
$$
We assume that $a, b \geq 0$. We put $p = a/(a+b)$ and $s^{1/(a+b)} = t$. Then 
$$
(t^{2a}-1)(t^{2b}-1) \geq \frac{4ab}{(a+b)^2}(t^{a+b}-1)^2.
$$
Then we have 
$$
(t^{2a}-1)(t^{2b}-1)(t^c+1)^2 \geq \frac{4ab}{(a+b)^2}(t^{a+b}-1)^2(t^c+1)^2. 
$$
In order to show the aimed inequality, we have to prove that 
$$
(t^{a+b}-1)^2(t^c+1)^2 \geq \frac{4(a+b)^2}{(a+b+c)^2}(t^{a+b+c}-1)^2.
$$
Since $a+b+c > 0$, it is sufficient to prove the following inequality
\begin{equation}
(t^{a+b}-1)(t^c+1) \geq \frac{2(a+b)}{a+b+c}(t^{a+b+c}-1)
\label{eq:num4-7}
\end{equation}
for $t \geq 1$ and $a, b, c \geq 0, 0 < a+b \leq c$ or $a, b \geq 0, c \leq 0, a+b+c > 0$. 
We put 
$$
S(t) = (t^{a+b}-1)(t^c+1)-\frac{2(a+b)}{a+b+c}(t^{a+b+c}-1).
$$
Then 
$$
S^{'}(t) = t^{c-1} \{(c-a-b)t^{a+b}-c+(a+b)t^{a+b-c} \}.
$$
Here we put 
$$
T(t) = (c-a-b)t^{a+b}-c+(a+b)t^{a+b-c}.
$$
Then 
$$
T^{'}(t) = (a+b)(c-a-b)t^{a+b-c-1}(t^c-1).
$$
When $a+b \leq c$, $T^{'}(t) \geq 0$. Since $T(1) = 0$, $T(t) \geq 0$ for $t \geq 1$. 
Then $S^{'}(t) \geq 0$. Since $S(1) = 0$, $S(t) \geq 0$ for $t \geq 1$.
On the other hand when $c \leq 0$, $T^{'}(t) \geq 0$.  Since $T(1) = 0$, $T(t) \geq 0$ for $t \geq 1$. 
Then $S^{'}(t) \geq 0$. Since $S(1) = 0$, $S(t) \geq 0$ for $t \geq 1$. Hence we get (\ref{eq:num4-7}). 
\ \hfill $\Box$ 

\begin{flushleft}
{\bf Proof of Proposition \ref{prop:proposition4-1}.}  Let $x < y$. In the last line of (\ref{eq:num4-5}), 
dividing both the numerator and the denominator by $(f(x)g(x)h(x))^2$ and by using 
$F(x) = \log f(x), G(x) = \log g(x)$ and $H(x) = \log h(x)$, we get 
\end{flushleft}
$$
L(x,y) = \frac{(e^{2(F(y)-F(x))}-1)(e^{2(G(y)-G(x))}-1)(e^{H(y)-H(x)}+1)^2}{(e^{F(y)-F(x)+G(y)-G(x)+H(y)-H(x)}-1)^2}  
$$
By the generalized mean value theorem, there exist $z~(x < z < y)$, $w~(x < w < y)$ such that 
$$
\frac{G(y)-G(x)}{F(y)-F(x)} = \frac{G^{'}(z)}{F^{'}(z)} = k(z),~~~\frac{H(y)-H(x)}{F(y)-F(x)} = \frac{H^{'}(w)}{F^{'}(w)} = \ell(w).
$$
Thus we have 
$$
L(x,y) = \frac{(e^{2(F(y)-F(x))}-1)(e^{2k(z)(F(y)-F(x))}-1)(e^{\ell(w)(F(y)-F(x))}+1)^2}{(e^{(1+k(z)+\ell(w))(F(y)-F(x))}-1)^2}. 
$$
It follows from Lemma \ref{lem:lemma4-1} that for any $R > 0$, the function 
$$
(k,\ell) \rightarrow A(k,\ell) = \frac{(R^2-1)(R^{2k}-1)(R^{\ell}+1)^2}{(R^{(1+k+\ell)}-1)^2}
$$
defined in $k \in [m,M], \ell \in [n,N]$ is bounded from below by $\min_{m \leq k \leq M, n \leq \ell \leq N} \{A(k,\ell) \}$. 
It is easy to obtain 
$$
\min_{m \leq k \leq M, n \leq \ell \leq N} \{A(k,\ell) \} \geq 16 \beta(f,g,h).
$$
We complete the proof. \ \hfill $\Box$

\begin{flushleft}
{\bf Proof of Theorem \ref{th:theorem3-1}.}  Since
\begin{eqnarray*}
Tr[f(\rho)g(\rho)h(\rho)[A,B]] & = & Tr[f(\rho)g(\rho)h(\rho)[A_0,B_0]] \\
& = & 2i {\rm Im}\{ Tr[f(\rho)g(\rho)h(\rho)A_0B_0]\} \\
& = & 2i {\rm Im} \sum_{\ell < m} (f(\lambda_{\ell})g(\lambda_{\ell})h(\lambda_{\ell})-f(\lambda_m)g(\lambda_m)h(\lambda_m))a_{m \ell}b_{\ell m} \\
& = & 2i \sum_{\ell < m} (f(\lambda_{\ell})g(\lambda_{\ell})h(\lambda_{\ell})-f(\lambda_m)g(\lambda_m)h(\lambda_m)){\rm Im}(a_{m \ell}b_{\ell m})
\end{eqnarray*}
for any $A, B \in M_{n,sa}(\mathbb{C})$, where $a_{\ell m} = \langle \phi_m|A_0|\phi_{\ell}\rangle$ and 
$b_{m \ell} = \langle \phi_m|B_0|\phi_m\rangle$, we have 
\begin{eqnarray*}
|Tr[f(\rho)g(\rho)h(\rho)[A,B]] & \leq & 2 \sum_{\ell < m} |f(\lambda_{\ell})g(\lambda_{\ell})h(\lambda_{\ell})-f(\lambda_m)g(\lambda_m)h(\lambda_m)| |{\rm Im}a_{m \ell} b_{m \ell}| \\
& \leq & 2 \sum_{\ell < m} |f(\lambda_{\ell})g(\lambda_{\ell})h(\lambda_{\ell})-f(\lambda_m)g(\lambda_m)h(\lambda_m)| |a_{m \ell}||b_{m \ell}|.
\end{eqnarray*}
By Proposition \ref{prop:proposition4-1}, we have 
\begin{eqnarray*}
&   & \beta(f,g,h)|Tr[f(\rho)g(\rho)h(\rho)[A,B]]|^2 \\
& \leq & 4 \beta(f,g,h) ( \sum_{\ell < m}|f(\lambda_{\ell})g(\lambda_{\ell})h(\lambda_{\ell})-f(\lambda_m)g(\lambda_m)h(\lambda_m)| |a_{m \ell}| |b_{\ell m}| )^2 \\
& \leq & \frac{1}{4} ( \sum_{\ell < m} \sqrt{(f(\lambda_{\ell})^2-f(\lambda_m)^2)(g(\lambda_{\ell})^2-g(\lambda_m)^2)(h(\lambda_{\ell})+h(\lambda_m))^2}|a_{\ell m}| |b_{m \ell}| )^2 \\
& = & \frac{1}{4} ( \sum_{\ell < m} \sqrt{ \Delta_f(\ell,m) \Delta_g(\ell,m) \Gamma_h(\ell,m)} |a_{m \ell}| \sqrt{ \Gamma_f(\ell,m) \Gamma_g(\ell,m) \Gamma_h(\ell,m)} |b_{\ell m}| )^2,
\end{eqnarray*}
where $\Delta_f(\ell,m) = f(\lambda_{\ell})-f(\lambda_m), \Delta_g(\ell,m) = g(\lambda_{\ell})-g(\lambda_m)$ and 
$\Gamma_f(\ell,m) = f(\lambda_{\ell})+f(\lambda_m), \Gamma_g(\ell,m) = g(\lambda_{\ell})+g(\lambda_m), 
\Gamma_h(\ell,m) = h(\lambda_{\ell})+h(\lambda_m)$.  By Schwarz inequality, we have 
\begin{eqnarray*}
&   & \beta(f,g,h)|Tr[f(\rho)g(\rho)h(\rho)[A,B]]|^2 \\
& \leq & \frac{1}{2} \sum_{\ell < m} \Delta_f(\ell,m) \Delta_g(\ell,m) \Gamma_h(\ell,m) |a_{m \ell}|^2 \\
&   & \times \frac{1}{2} \sum_{\ell < m} \Gamma_f(\ell,m) \Gamma_g(\ell,m) \Gamma_h(\ell,m) |b_{\ell m}|^2 \\
& \leq & I_{\rho,(f,g,h)}(A) J_{\rho,(f,g,h)}(B).
\end{eqnarray*}
Similarly we have 
$$
\beta(f,g,h)|Tr[f(\rho)g(\rho)h(\rho)[A,B]]|^2 \leq I_{\rho,(f,g,h)}(B) J_{\rho,(f,g,h)}(A).
$$
Hence by multiplying the above two inequalities, we have 
$$
\beta(f,g,h)|Tr[f(\rho)g(\rho)h(\rho)[A,B]]|^2 \leq U_{\rho,(f,g,h)}(A) U_{\rho,(f,g,h)}(B).
$$
\end{flushleft}
\ \hfill $\Box$

When $h(x) = 1$, we obtain the result given by Ko and Yoo \cite{KoYo:monotone}. 

\begin{corollary}[\cite{KoYo:monotone}]
If $(f,g)$ is CLI monotone pair, then the following inequality holds: 
$$
U_{\rho,(f,g)}(A) U_{\rho,(f,g)}(B) \geq \beta(f,g)|Tr[f(\rho)g(\rho)[A,B]]|^2 
$$
for $A, B \in M_{n,sa}(\mathbb{C})$, where 
\begin{eqnarray*}
I_{\rho,(f,g)}(A) & = & \frac{1}{2}Tr[(i[f(\rho),A_0])(i[g(\rho),A_0])], \\
J_{\rho,(f,g)}(A) & = & \frac{1}{2}Tr[\{ f(\rho),A_0 \} \{ g(\rho),A_0 \}], \\
U_{\rho,(f,g)}(A)& = & \sqrt{ I_{\rho,(f,g)} J_{\rho,(f,g)}}, \\
\beta(f,g) & = & \min \{ \frac{m}{(m+M)^2},~ \frac{M}{(m+M)^2} \}.
\end{eqnarray*}
\label{cor:corollary4-1}
\end{corollary}

We also have the following corollary.         

\begin{corollary}
Let $f(x) = x^{\alpha}~(\alpha \geq 0),~~g(x) = x^{\beta}~(\beta \geq 0),~~h(x) = x^{\gamma}~(\gamma \geq 0~{\rm or}~\gamma \leq 0)$.
\begin{description}
\item[(1)]  If $\alpha, \beta, \gamma \geq 0$ satisfy $ 0 < \alpha + \beta \leq \gamma$, then 
$$
\beta(f,g,h) = \frac{\alpha \beta}{(\alpha+\beta+\gamma)^2}.
$$
\item[(2)] If $\alpha, \beta \geq 0, \gamma \leq 0$ satisfy $\alpha+\beta+\gamma > 0$, then 
$$
\beta(f,g,h) = \frac{\alpha \beta}{(\alpha+\beta+\gamma)^2}.
$$
\end{description} 
\label{cor:corollary4-2}
\end{corollary}

\begin{remark}
When $\alpha, \beta \geq 0, \gamma < 0$ satisfy $\alpha + \beta + \gamma > 0$,  we remark that 
$h(x)$ is not continuous function on $[0,1]$ because
$$
\lim_{x \to +0} h(x) = +\infty.
$$
Then in this case by putting $\epsilon > 0$ such that $\epsilon$ is smaller than the minimal eigenvalue of $\rho$,  
we can assume that $h(x)$ is continuous on $[\epsilon,1]$. Hence we obtain the same result as Corollary \ref{cor:corollary4-2}. 
\label{rem:remark4-1}
\end{remark}

\begin{remark}
When $\gamma = 0$ in (2) of Corollary \ref{cor:corollary4-2}, we have the result in \cite {Ya:GWYD2}~(Theorem \ref{th:theorem2-3}). 
And when $\alpha+\beta+\gamma = 1$ in Corollary \ref{cor:corollary4-2}, we have the result in \cite{Ya:GWYD1}~(Theorem \ref{th:theorem2-2}). 
That is (1) implies $\alpha, \beta \geq 0, \alpha+\beta \leq \frac{1}{2}$ and (2) implies $\alpha, \beta \geq 0, \alpha+\beta \geq 1$.
\label{rem:remark4-2}
\end{remark}

\end{document}